\documentclass[runningheads]{llncs}
\usepackage[breaklinks, colorlinks, linkcolor=blue, citecolor=blue, urlcolor=black]{hyperref}
\usepackage{color}

\usepackage{graphicx}
\usepackage{tabularx}
\usepackage{tabularcalc}
\usepackage{booktabs}
\usepackage{multirow}
\usepackage{algorithm}
\usepackage{algpseudocode}
\usepackage{tikz}
\usepackage{forest}
\usepackage{pgfplots}
\usepackage{xltabular}
\usetikzlibrary{arrows.meta}
\usepackage{amsmath}
\usepackage{amssymb}
\usepackage{xspace}
\usepackage{xcolor}
\usepackage{varioref}
\usepackage{aligned-overset}
\usepackage[shortlabels,inline]{enumitem}
\usepackage{boxedminipage}
\usepackage{makecell}
\usepackage{underscore}
\usepackage{rotating}
\usepackage{cleveref}
\usepackage{hyperref}

\usepackage{pgfplots}
\usepackage{booktabs}
\pgfplotsset{compat=1.18}

\definecolor{pbblue}{RGB}{68,119,170}
\definecolor{pbcyan}{RGB}{102,204,238}
\definecolor{pbgreen}{RGB}{34,136,51}
\definecolor{pbred}{RGB}{238,102,119}

\usepackage{ifthen}

\usepackage{ifthen}

\provideboolean{detectedSTOC}
\provideboolean{detectedFOCS}
\provideboolean{detectedElsevier}
\provideboolean{detectedNOW}
\provideboolean{detectedLMCS}
\provideboolean{detectedIEEE}
\provideboolean{detectedPoster}
\provideboolean{detectedSIAM}
\provideboolean{detectedLNCS}
\provideboolean{detectedACM}
\provideboolean{detectedACMconf}
\provideboolean{detectedSigplanconf}
\provideboolean{detectedToC}
\provideboolean{detectedLIPIcs}
\provideboolean{detectedAAAI}
\provideboolean{detectedIJCAI}
\provideboolean{detectedCompCplx}
\provideboolean{detectedEasyChair}
\provideboolean{detectedArticle}
\provideboolean{detectedReport}
\provideboolean{detectedThesis}

\makeatletter

\@ifclassloaded{sig-alternate}
{\setboolean{detectedSTOC}{true}}
{\setboolean{detectedSTOC}{false}}

\@ifclassloaded{elsarticle}
{\setboolean{detectedElsevier}{true}}
{\setboolean{detectedElsevier}{false}}

\@ifclassloaded{now}
{\setboolean{detectedNOW}{true}}
{\setboolean{detectedNOW}{false}}

\@ifclassloaded{lmcs}
{\setboolean{detectedLMCS}{true}}
{\setboolean{detectedLMCS}{false}}

\@ifclassloaded{IEEEtran} {
  \setboolean{detectedIEEE}{true}
  \ifCLASSOPTIONconference {
    \setboolean{detectedFOCS}{true}
  }
  \else {
    \setboolean{detectedFOCS}{false}
  }
  \fi
}
{
  \setboolean{detectedFOCS}{false}
  \setboolean{detectedIEEE}{false}
}

\@ifclassloaded{siamart171218}
{\setboolean{detectedSIAM}{true}}
{\setboolean{detectedSIAM}{false}}

\@ifclassloaded{llncs}
{\setboolean{detectedLNCS}{true}}
{\setboolean{detectedLNCS}{false}}

\@ifclassloaded{acmsmall}
{\setboolean{detectedACM}{true}}
{\setboolean{detectedACM}{false}}

\@ifclassloaded{acmart}
{\setboolean{detectedACMconf}{true}}
{\setboolean{detectedACMconf}{false}}

\@ifclassloaded{sigplanconf}
{\setboolean{detectedSigplanconf}{true}}
{\setboolean{detectedSigplanconf}{false}}

\@ifclassloaded{toc}
{\setboolean{detectedToC}{true}}
{\setboolean{detectedToC}{false}}

\@ifclassloaded{lipics}
{\setboolean{detectedLIPIcs}{true}}
{\@ifclassloaded{lipics-v2019}
  {\setboolean{detectedLIPIcs}{true}}
  {\@ifclassloaded{oasics-v2019}
    {\setboolean{detectedLIPIcs}{true}}
    {\setboolean{detectedLIPIcs}{false}}
  }
}

\@ifclassloaded{cc}
{\setboolean{detectedCompCplx}{true}}
{\setboolean{detectedCompCplx}{false}}

\@ifclassloaded{easychair}
{\setboolean{detectedEasyChair}{true}}
{\setboolean{detectedEasyChair}{false}}

\@ifpackageloaded{aaai}
{\setboolean{detectedAAAI}{true}}        
{\@ifpackageloaded{aaai18}
  {\setboolean{detectedAAAI}{true}}
  {\@ifpackageloaded{aaai20}       
    {\setboolean{detectedAAAI}{true}}
    {\setboolean{detectedAAAI}{false}}}}

\@ifpackageloaded{ijcai18}
{\setboolean{detectedIJCAI}{true}}        
{\@ifpackageloaded{ijcai19}
  {\setboolean{detectedIJCAI}{true}}        
  {\setboolean{detectedIJCAI}{false}}}

\@ifclassloaded{sciposter}
{\setboolean{detectedPoster}{true}}
{\setboolean{detectedPoster}{false}}

\@ifclassloaded{article}
{\setboolean{detectedArticle}{true}}
{\setboolean{detectedArticle}{false}}

\makeatother

\ifthenelse{\not \isundefined{\examen} 
  \and \not \isundefined{\disputationsdatum} 
  \and \not \isundefined{\disputationslokal}}   
  {\setboolean{detectedThesis}{true}}
  {\setboolean{detectedThesis}{false}}

\ifthenelse{\boolean{detectedArticle} \or \boolean{detectedThesis}
  \or \boolean{detectedSTOC} \or \boolean{detectedFOCS}
  \or \boolean{detectedSIAM} \or \boolean{detectedIEEE}
  \or \boolean{detectedPoster}}
{\setboolean{detectedReport}{false}}
{\setboolean{detectedReport}{true}}

\usepackage{xspace}
\ifthenelse
{\boolean{detectedElsevier} \or \boolean{detectedSIAM} 
  \or \boolean{detectedLIPIcs}}
{}
{\usepackage{varioref}}

\makeatletter
\AtBeginDocument{%
  \def\doi#1{\url{https://doi.org/#1}}}
\makeatother

\DeclareMathAlphabet{\mathsfsl}{OT1}{cmss}{m}{sl}

\DeclareRobustCommand{\BibTeX}{%
  {\normalfont B\kern-.05em{\scshape i\kern-.025em b}\kern-.08em \TeX}%
}

\ifthenelse{\boolean{detectedToC}}{}
{

}

\newcommand{\MAXOFEXPR}[2][]{\max_{#1} \left\{ #2 \right\}}
\newcommand{\MINOFEXPR}[2][]{\min_{#1} \left\{ #2 \right\}}
\newcommand{\Maxofexpr}[2][]{\max_{#1} \bigl\{ #2 \bigr\}}
\newcommand{\Minofexpr}[2][]{\min_{#1} \bigl\{ #2 \bigr\}}

\newcommand{\MAXOFSET}[3][:]%
     {\ifthenelse{\equal{#1}{;}}%
     {\MAXOFEXPR{ #2 \,;\, #3 }}
     {\ifthenelse{\equal{#1}{:}}%
     {\MAXOFEXPR{ #2 \,:\, #3 }}
     {\max \twincommandJN{\left\{}{#2}{\left#1}{\right}{\,#3}{\right\}}}}}
\newcommand{\MINOFSET}[3][:]%
     {\ifthenelse{\equal{#1}{;}}%
     {\MINOFEXPR{ #2 \,;\, #3 }}
     {\ifthenelse{\equal{#1}{:}}%
     {\MINOFEXPR{ #2 \,:\, #3 }}
     {\min \twincommandJN{\left\{}{#2}{\left#1}{\right}{\,#3}{\right\}}}}}

\newcommand{\Maxofset}[3][:]%
     {\ifthenelse{\equal{#1}{;}}%
     {\Maxofexpr{ #2 \,;\, #3 }}
     {\ifthenelse{\equal{#1}{:}}%
     {\Maxofexpr{ #2 \,:\, #3 }}
     {\max \twincommandJN{\bigl\{}{#2}{\bigl#1}{\bigr}{\,#3}{\bigr\}}}}}
\newcommand{\Minofset}[3][:]%
     {\ifthenelse{\equal{#1}{;}}%
     {\Minofexpr{ #2 \,;\, #3 }}
     {\ifthenelse{\equal{#1}{:}}%
     {\Minofexpr{ #2 \,:\, #3 }}
     {\min \twincommandJN{\bigl\{}{#2}{\bigl#1}{\bigr}{\,#3}{\bigr\}}}}}

\DeclareMathOperator{\Expop}{E}

\ifthenelse{\boolean{detectedLMCS}}
{}
{}

\newcommand{\twincommandJN}[6]%
    {#1#2#3\vphantom{#2#5}\mspace{-2.05mu}#4.#5#6}

\newcommand{\CondExp}[2]%
    {\Expop\twincommandJN{\bigl[}{#1}{\bigl|}{\bigr}{\,#2}{\bigr]}}
\newcommand{\CONDEXP}[2]%
     {\Expop\twincommandJN{\left[}{#1}{\left|}{\right}{\,#2}{\right]}}

\newcommand{\Condprob}[3][]%
    {\Pr_{#1}\twincommandJN{\bigl[}{#2}{\bigl|}{\bigr}{\,#3}{\bigr]}}
\newcommand{\CONDPROB}[3][]%
    {\Pr_{#1}\twincommandJN{\left[}{#2}{\left|}{\right}{\,#3}{\right]}}

\newcommand{\Set}[1]{\bigl\{ #1 \bigr\}}

\newcommand{\Setdescr}[3][|]%
     {\ifthenelse{\equal{#1}{;}}%
     {\Set{ #2 \,;\, #3 }}
     {\ifthenelse{\equal{#1}{:}}%
     {\Set{ #2 \,:\, #3 }}
     {\twincommandJN{\bigl\{}{#2\,}{\bigl#1}{\bigr}{\,#3}{\bigr\}}}}}
\newcommand{\SETDESCR}[3][|]%
     {\twincommandJN{\left\{}{#2\,}{\left#1}{\right}{\,#3}{\right\}}}

\newcommand{\Setdescrbrackets}[3][|]%
     {\twincommandJN{\bigl[}{#2}{\bigl#1}{\bigr}{\,#3}{\bigr]}}
\newcommand{\SETDESCRBRACKETS}[3][|]%
     {\twincommandJN{\left[}{#2}{\left#1}{\right}{\,#3}{\right]}}

\newcommand{\nvar}{n}

\newcommand{\nclause}{m}

\newcommand{\clwidth}{k}

\newcommand{\randkcnfnclwrepl}[3][\clwidth]%
        {\ensuremath{\mathcal{F}^{#2, #3}_{#1}}}
\newcommand{\randkcnfnclwreplstd}%
        {\randkcnfnclwrepl{\clwidth}{\nvar}{\nclause}}

\newcommand{\complclassformat}[1]%
        {\textrm{\upshape{\textsf{#1}}}\xspace}
\newcommand{\cocomplclass}[1]%
        {\textrm{\upshape{\textsf{co#1}}}\xspace}

\newcommand{\DTIMEadviceclass}[2]%
    {\ensuremath{\complclassformat{DTIME}\bigl(#1\bigr)/{#2}}}

\newcommand{\PCPalph}[5]%
    {\ensuremath{\complclassformat{PCP}_{{#1},{#2}}[{#3}, {#4}, {#5}]}}
\newcommand{\PCP}[4]%
    {\ensuremath{\complclassformat{PCP}_{{#1},{#2}}[{#3}, {#4}]}}

\ifthenelse{\boolean{detectedToC}}{}
  { 
    }
\ifthenelse{\boolean{detectedToC}}{}{

}

\ifthenelse{\boolean{detectedLIPIcs} \or \boolean{detectedIJCAI}}
{\renewcommand{\st}{\errmessage{Please do not use st}}}
{\ifthenelse{\isundefined{\st}}
  {\newcommand{\st}{such that\xspace}}
  {}}      

\ifthenelse{\boolean{detectedIEEE}}{}{}

\ifthenelse
{\isundefined{\refeq}}
{\newcommand{\refeq}[1]{\eqref{#1}}}
{\renewcommand{\refeq}[1]{\eqref{#1}}}

\ifthenelse{\boolean{detectedACMconf}}
{}
{}

\ifthenelse{\boolean{detectedACMconf}}
{
 }
{
 }

\newcommand{\SETSOFVARSORLIT}[2]%
        {\mathit{#1}\left({#2}\right)}
\newcommand{\setsofvarsorlit}[2]%
        {\mathit{#1}({#2})}
\newcommand{\Setsofvarsorlit}[2]%
        {\mathit{#1}\bigl({#2}\bigr)}

\newcommand{\derivabbrev}[2]{\bigl( #1 \vdash #2 \bigr)}
\newcommand{\derivabbrevsmall}[2]{( #1 \vdash #2 )}
\newcommand{\derivabbrevcompact}[2]{\bigl( #1 \vdash #2 \bigr)}

\newcommand{\refutabbrevsmall}[1]{\derivabbrevsmall{#1}{\!\bot}}
\newcommand{\refutabbrevcompact}[1]{\derivabbrevcompact{#1}{\!\bot}}

\newcommand{\genericrefsmall}[3]%
    {{\mathit{#1}}_{#2}\refutabbrevsmall{#3}}
\newcommand{\genericrefcompact}[3]%
    {{\mathit{#1}}_{#2}\refutabbrevcompact{#3}}
\newcommand{\genericderiv}[4]%
    {{\mathit{#1}}_{#2}\derivabbrev{#3}{#4}}
\newcommand{\genericderivsmall}[4]%
    {{\mathit{#1}}_{#2}\derivabbrevsmall{#3}{#4}}
\newcommand{\genericderivcompact}[4]%
    {{\mathit{#1}}_{#2}\derivabbrevcompact{#3}{#4}}
\newcommand{\generictaut}[3]%
    {{\mathit{#1}}_{#2}\derivabbrev{}{#3}}
\newcommand{\generictautcompact}[3]%
    {{\mathit{#1}}_{#2}\derivabbrevcompact{}{#3}}
\newcommand{\generictautsmall}[3]%
    {{\mathit{#1}}_{#2}\derivabbrevsmall{}{#3}}

\newcommand{\formulaformat}[1]{\mathit{#1}}

\newcommand{\extendedversion}[1]{\widetilde{#1}}

\newcommand{\epopnot}[1]%
    {\extendedversion{\formulaformat{POP}}_{#1}}

\newcommand{\elopnot}[1]%
    {\extendedversion{\formulaformat{LOP}}_{#1}}

\newcommand{\ephpnot}[2]%
    {\vphantom{\extendedversion{\formulaformat{PHP}}}
      {\smash{\extendedversion{\formulaformat{PHP}}}
        \vphantom{\formulaformat{PHP}}}^{#1}_{#2}}

\newcommand{\efphpnot}[2]%
    {\vphantom{\extendedversion{\formulaformat{FPHP}}}
      {\smash{\extendedversion{\formulaformat{FPHP}}}
        \vphantom{\formulaformat{FPHP}}}^{#1}_{#2}}

\newcommand{\ontophpnot}[2]%
    {\formulaformat{Onto}\text{-}\formulaformat{PHP}^{#1}_{#2}}
\newcommand{\ontofphpnot}[2]%
    {\formulaformat{Onto}\text{-}\formulaformat{FPHP}^{#1}_{#2}}

\newcommand{\graphontophpnot}[1][G]%
    {\text{$\formulaformat{Onto}$-$\formulaformat{PHP}$}({#1})}

\newcommand{\perfectmatchingnot}[1][G]%
    {\formulaformat{PM}({#1})}
\newcommand{\RNum}[1]{\uppercase\expandafter{\romannumeral #1\relax}}
\newcommand{\solver}[1]{\textsc{#1}\xspace}

\newcommand{\scip}{\solver{SCIP}}
\newcommand{\highs}{\solver{HiGHS}}

\newcommand{\papilo}{\solver{PaPILO}}

\newcommand{\miplib}{\solver{MIPLIB}}

\newcommand{\myorcidlink}[1]{\,\href{https://orcid.org/#1}{\raisebox{-0.45ex}{\includegraphics[width=1.8ex]{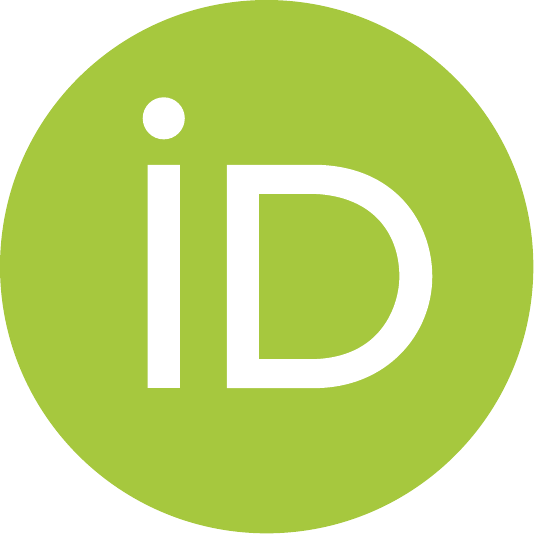}}}}

\definecolor{pbblue}{RGB}{68,119,170}
\definecolor{pbcyan}{RGB}{102,204,238}
\definecolor{pbgreen}{RGB}{34,136,51}
\definecolor{pbred}{RGB}{238,102,119}

\usepackage{todonotes,xspace}

\makeatletter
\def\orcidID#1{\href{http://orcid.org/#1}{\protect\raisebox{-1.25pt}{\protect\includegraphics{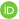}}}}
\makeatother

\makeatletter
\@ifpackageloaded{hyperref}{%
  \newcounter{algorithmicH}%
  \setcounter{algorithmicH}{0}%
  \renewcommand{\theHALG@line}{%
    \thealgorithm.\arabic{ALG@line}%
  }%
}{}
\makeatother

\begin{document}

\title{Clique Probing for Mixed-Integer Programs}
\author{
    Jacob von Holly-Ponientzietz\inst{1}\orcidID{0009-0002-2601-3689} \and
    Alexander Hoen\inst{1,2}\orcidID{0000-0003-1065-1651} \and
    Mark Turner\inst{1}\orcidID{0000-0001-7270-1496} \and
    Ambros Gleixner\inst{1,2}\orcidID{0000-0003-0391-5903}
    \institute{%
      Zuse Institute Berlin, Takustr.~7, 14195 Berlin, Germany\\
      \email{von.holly-ponientzietz@zib.de,turner@zib.de}
      \and
      HTW Berlin, 10313 Berlin, Germany\\
      \email{hoen@htw-berlin.de,gleixner@htw-berlin.de} \\
    }
    }

\authorrunning{von Holly et. al}
\titlerunning{Clique probing} 

\maketitle

\begin{abstract}
  Probing is an important presolving technique in mixed-integer programming solvers.
It selects binary variables, tentatively fixes them to 0 and 1, and performs propagation to deduce additional variable fixings, bound tightenings, substitutions, and implications. 
In this work, we propose \emph{clique probing}: instead of probing on individual variables, we select cliques, a set of binary variables of which at most one can be set to one, and systematically probe on all variables of a clique.
Experiments with our implementation in the open-source presolve library \papilo demonstrate that exploiting clique information in this form significantly increases the number of reductions.
When integrated into the MIP solver \scip, we observe a 3\% performance improvement on \miplib instances containing cliques.
\end{abstract}

\section{Introduction}
\label{sec::introduction}

\emph{Mixed Integer Programming} (MIP) is widely used in industry to solve optimization problems, with competitive MIP solvers 
typically implementing some variant of \textit{branch-\&-cut} 
as the backbone of the solving process.
Before initializing the branch-\&-cut process, however, \emph{presolving} (or \emph{preprocessing}) is called.
The goal of presolving is to eliminate redundant information and strengthen the problem formulation, 
thereby accelerating the subsequent solution process. 
Therefore, presolving is a crucial part of modern MIP solvers and, more often than not, it is the deciding factor between a problem being solvable or not~\cite{presolving_achterberg,achterberg2013mixed}. 


One important presolving technique is \emph{probing}~\cite{Savelsbergh_Preprocessing_and_Probing_Techniques,Achterberg07Thesis}.
The central idea of probing is to select binary variables and temporarily fix them to 0 and 1. 
After each tentative fixing, a propagation step is performed, and the resulting deductions are analyzed in order to strengthen the formulation.
In \cite{presolving_achterberg}, disabling probing led to a slowdown of 7\%  to 33\% on instances with more than 1000~seconds solving time and rendered 11 of 3155 models unsolvable.

\begin{definition}
\label{def::probing}
Formally, let $x$ be a binary variable and $z$ an arbitrary variable with bounds $[l_z,u_z]$.
We denote the variable bounds of $z$ after propagating $x=0$ as $[l_z^0,u_z^0]$ and the variable bounds of $z$ after propagating $x=1$ as $[l_z^1,u_z^1]$.
The following four rules can be applied while probing:
\begin{description} 
    \item[inf] If setting $x = 0$ ($x=1$) makes the problem infeasible, $x$ can be fixed to 1 (0).
    \item[bds] Global bounds for $z$ can be updated to $[\min(l_z^0,l_z^1), \max(u_z^0, u_z^1)]$.
    \item[sub] If $l^0_z = u^0_z$ and $l^1_z = u^1_z$, the variable $z$ can be substituted by $z = l^0_z + (l^1_z - l^0_z)\; x$.
    \item[imp] Implications can be collected in a central data structure of the solver. An example for an implication would be $x_k = 0 \rightarrow z \geq l_z^0$. 
\end{description}
\end{definition}

Performing probing on all variables can be time-consuming and may outweigh its potential benefits to the solver.
To address this, probing orders the variables according to a scoring function and employs a termination criterion to stop early.
Thus, designing an effective scoring function and appropriate termination criteria is crucial for an efficient implementation.
Despite the importance of both probing and the subsequent solving process, we are not aware of any published work describing scoring methods or termination criteria.

Open-source solvers like \scip~\cite{SCIP9} and \highs~\cite{HH2015parallelizingdualrevisedsimplex} assign scores to each binary variable and rank them by this score. The score tries to reflect the impact that fixing the variables has on the model.
These solvers abort if a certain absolute or relative limit is reached, typically a limit on time, number of variables, number of fixings, or a combination of these limits.
In order to prevent stalling, an additional termination criterion is typically added, which interrupts probing if the last iterations were barely successful. 
MIP solvers typically exploit problem structure to accelerate the solving process.
Examples include the use of clique tables~\cite{Nemhauser1975VertexPacking,Achterberg07Thesis}, specialized handling of SOS constraints, and the identification of knapsack structure for generating strong cover cuts~\cite{Balas1975KnapsackFacets,Atamturk2004CoverInequalities,Achterberg07Thesis}.
However, to the best of our knowledge, this has not yet been used for probing. 
To improve the efficiency of probing, we propose leveraging structural information from the model. In particular, we focus on cliques to speed up the procedure, although the underlying idea can be generalized to other types of structure.
Therefore, we propose \emph{Clique Probing}: 
During our probing routine, we use the information provided by cliques to adjust the fixings done during probing. 
By this, we get stronger reductions and are able to speed up the subsequent solving process.

The paper is organized as follows: 
In Section \ref{sec::cliqueprobing}, we present the idea of clique probing 
and its implementation.
In Section \ref{sec::experiments}, we evaluate the performance of clique probing, and summarize our results in Section \ref{sec::conclusion}.

\section{Clique Probing}
\label{sec::cliqueprobing}

\subsection{General idea}


A \emph{clique} is a set of $n$ binary variables $x_1,\ldots,x_n$ of which at most one variable is allowed to be active, i.e., to be set to one.  We distinguish \emph{at-most-one cliques}, for which
\begin{align*}
    x_1 + \dots + x_n \leq 1
\end{align*}
must hold, and \emph{exactly-one cliques}, for which
\begin{align*}
    x_1 + \dots + x_n  = 1
\end{align*}
must hold.
During presolving, cliques are found explicitly as part of the constraints of the original model or derived from other problem structures that imply mutually exclusive assignments.

Standard probing selects variables based on the score, fixes them to one and zero, and performs propagation on both assignments.
Assuming all variables of a clique are selected for probing, probing performs $2n$ propagation calls, two for each variable.
However, a clique only allows at most $n+1$ different assignments for the variable vector $(x_1,\ldots,x_n)$, and for exactly-one cliques, even one less, because the all-zero assignment is excluded.
Therefore, it would be sufficient to just probe these $n+1$, respectively $n$, cases.
This has the advantage that we can
\begin{enumerate*}[label=(\alph*)]
    \item reduce the number of probing and subsequent propagation calls needed to get the same results,
    \item propagate faster since we can fix multiple variables initially,
    \item obtain tighter bounds for cliques involving more than two variables, since fixing a variable inside the clique to zero does not determine the remaining variables within the clique, and
    \item exploit that the all-zero assignment is not feasible for exactly-one cliques.
\end{enumerate*}

The following definition adjusts the evaluation of the results of standard probing in Definition~\ref{def::probing} for probing on cliques.


\begin{definition}
    Formally, let $c$ be an at-most-one clique consisting of $n$ variables.
Let $c^i$ be the variable assignment, where $x_i$ is set to one, and all other variables of clique $c$ are set to zero.
Let $c^0$ be the variable assignment, where all variables of a clique $c$ are set to $0$.
Let $z$ be an arbitrary variable with bounds $[l_z,u_z]$, and let $[l_z^i,u_z^i]$ be the variable bounds of $z$ after propagating $c^i$.
Then a new set of four rules can be applied while probing:

\begin{description} 
    \item[inf] If propagating $c^i$ with $i \in \{1,\dots,n\}$ makes the problem infeasible, $x_i$ can be fixed to zero. If propagating $c^0$ makes the problem infeasible, $c$ can be upgraded to an exactly-one clique.
    \item[bds] Global bounds for $z$ can be updated to $[\min_{k \in \{0,\dots,n\}}\;\ell^k_z, \max_{k \in \{0,\dots,n\}}\;u^k_z]$.
    \item[sub] Variable $z$ can be substituted by $z = \ell_z^j + ( \ell_z^i - \ell^j_z)\; x_i$ with $j \in \{1,\dots,n\}\backslash\{i\}$ if the following rules hold:
        \begin{enumerate}
            \item $\ell^k_z = u^k_z$ for all $ k\in \{0,1,\dots,n\}$
            \item $\ell^j_z = \ell^k_z$ and  $u^j_z = u^k_z$ for all $j,k\in \{0,1\dots,n\} \backslash\{i\}$
            \item $\ell^k_z \not= \ell^i_z$ and  $u^k_z \not= u^i_z$ for all $k\in\{0,1\dots,n\} \backslash\{i\}$ 
        \end{enumerate}
    \item[imp] Implications for $c^i$ with $i \in \{1,\dots,n\}$ can still be stored similarly to Definition \ref{def::probing}.
\end{description}
\end{definition}


For evaluating exactly-one cliques, the variable assignment $c^0$ is not considered, and the rules are adapted accordingly: in the \emph{bds} and \emph{subs} rules, the index sets range over $\{1,\dots,n\}$ instead of $\{0,\dots,n\}$. 
Further, the sub-rule can be more generalized, allowing the aggregation of multiple variables.

\subsection{Data structure}
\label{sec::data_structure}

In formulating the rules for analyzing substitutions and bounds, it is necessary to know all bounds of each variable after the propagation of variable assignments.
For long cliques, however, this can be memory-intensive and cause a significant slowdown.

To address this, we maintain a global data structure within clique probing instead of storing all bounds for each propagation.
For the global bounds, we record for each variable the current maximum upper and minimum lower bound, denoted as $B$ in Algorithm \ref{alg::probeclique} below.
To detect substitutions, we additionally store the second lowest and highest values, as well as the indices (argmin/argmax) of the minimum and maximum bounds, denoted as $I$ in Algorithm \ref{alg::probeclique}.

Rather than duplicating bounds for every possible variable assignment, this approach limits the overhead to six values per variable.
Additionally, after propagating each assignment, these values allow us to calculate the maximal bound changes and substitutions as if no further assignments were to be considered.
For example, if a variable’s bound is not tightened after propagating a single assignment, then no global tightening of that bound can occur.
These values can be used as an upper bound on the number of potential bound changes and substitutions, which can then be used as a termination criterion for large cliques.

\subsection{Algorithm}

Algorithm \ref{alg::probeclique} shows the function to perform probing on a single clique. 
First, the data structures $B$ (for tracking global bound changes) and $I$ (for detecting substitutions) are initialized (Line \ref{line::func::init}).
For at-most-one cliques, all variables in the clique are fixed to zero and propagated.
If this all-zero propagation is infeasible, the clique is upgraded to an exactly-one clique (Line \ref{line::func::upgrade}).
Otherwise, $B$ and $I$ are updated according to rules \emph{(bds, subs)} (Line \ref{line::func::update_zero}). 
Next, each variable in the clique is probed by setting it to one while fixing all others to zero, followed by propagation.
If infeasibility is detected, the corresponding variable is fixed to zero (Line \ref{line::func::fix_var}); otherwise, the data structures $B$ and $I$ are updated (Line \ref{line::func::update_vars}).

Since this process can become computationally expensive for large cliques, we introduce an early termination criterion.
As explained in Section \ref{sec::data_structure}, we can derive an upper bound on how many bound changes can be applied from $B$ and $I$.
If this falls below a certain threshold, we abort clique probing (Line \ref{line::func::check}).
Finally, after probing all cases, the results stored in $B$ and $I$ are analyzed, and the corresponding bound changes and substitutions are applied (Lines \ref{line::func::analyze_bounds} and \ref{line::func::analyze_subs}). 
If each propagation within the clique proves infeasible, the algorithm concludes global infeasibility (Line \ref{line::func::global_inf}).

\algtext*{EndIf}
\algtext*{EndFor}
\algtext*{EndWhile}
\algtext*{EndFunction}
\algdef{SE}[PARFOR]{ParFor}{EndParFor}[1]{\textbf{in parallel for} #1}{}
\algtext*{EndParFor}

\begin{algorithm}
    \caption{\textsc{ProbeSingleClique}}
    \label{alg::probeclique}
    \begin{algorithmic}[1]
            \Require MIP with global variable bound vectors $l,u$; clique $c$ with the set of all variables in the clique $N_c$; abort threshold $\delta$
            \State init $B,I$ with global bounds\label{line::func::init}
            \If{$c$ is an at-most-one clique}
                \State set all $x \in N_c$ to 0; propagate
                \If{infeasible} 
                    \State upgrade clique to an exactly-one clique\label{line::func::upgrade}
                \Else 
                    \State update $B, I$\label{line::func::update_zero}
                \EndIf
            \EndIf
            \For{$x \in N_c$}\label{line::func::for}
                \State set $x$ to 1; set all $y \in N_c/ \{x\}$ to 0; propagate
                \If{\Call{MaxGlobalBoundChangesAndImplications}{$B,I$} $\leq \delta$}
                    \label{line::func::check}
                    \State \textbf{abort}
                \EndIf
                \If{infeasible} 
                    \State fix $x$ to $0$\label{line::func::fix_var}
                \Else 
                    \State update $B, I$\label{line::func::update_vars}
                \EndIf
            \EndFor
            \State analyze $B$ and apply global valid variable bounds\label{line::func::analyze_bounds}
            \State analyze $I$ and apply substitutions\label{line::func::analyze_subs}
            \If{all propagations infeasible}
                \State \Return global infeasibility 
            \EndIf\label{line::func::global_inf}
    \end{algorithmic}
\end{algorithm}


A full run of clique probing is outlined in Algorithm \ref{alg::general_algorithm}.
First, we need to select the cliques to be probed on. 
In Line \ref{line::alg::scoring}, we sort cliques by their average variable score, and the highest-scoring cliques
Then, we perform clique probing on the selected cliques until either all the cliques have been probed or two consecutive cliques have been unsuccessful (Line \ref{line::alg::abortion}).
Probing on a clique is considered unsuccessful if a certain threshold of bound changes and substitutions per propagated variable is not reached.
Afterwards, traditional probing is continued in Line \ref{line::alg::probing} on the variables ($N_P$) that have not been touched by clique probing (Line \ref{line::alg::exclude_vars}).
Clique probing is disabled for future runs if it finds no reductions at all (Line \ref{line::alg::disable}).

\begin{algorithm}
    \caption{\textsc{CliqueProbing}}
    \label{alg::general_algorithm}
    \begin{algorithmic}[1]
            \Require a list of cliques $C$, probing score for binary variables $S$, list of binary variables in the MIP $N$
            \State Sort $C$ by average probing score of variables $S$ \label{line::alg::scoring}
            \State $C\gets $ highest-scored cliques, skipping cliques with more than 50\% overlapping variables \label{line::alg::select_cliques}
            \State $N_P \gets N$ 
            \For{$c \in \mathit{C}$}\label{line::alg::for}
                \State \Call{ProbeSingleClique}{$c$}
                \State $N_P \gets N_P \backslash \{N_c\}$ \label{line::alg::exclude_vars}
                \If{two consecutive calls of \Call{ProbeSingleClique}{} were unsuccessful} \label{line::alg::abortion}
                    \State \textbf{break}
                \EndIf
            \EndFor
            \State disable clique probing for future runs if no reductions were found \label{line::alg::disable}
            \State resume with \Call{Probing}{$N_P$}\label{line::alg::probing}
    \end{algorithmic}
\end{algorithm}

\subsection{Implementation}
\label{sec::implementation}

For our implementation, we build upon the open-source presolving C++ library \papilo~\cite{GGHpapilo}, which is publicly available on GitHub\footnote{\url{https://github.com/scipopt/papilo}}.
\papilo implements common presolving techniques~\cite{presolving_achterberg}, supports multi-precision, and accelerates the presolving process through parallelization.
To minimize synchronization overhead, it adopts a transaction-based architecture: each presolver records its reductions independently in the form of a transaction rather than applying them immediately. 
These transactions are then returned to the \papilo core, which sequentially checks their validity and applies them if still applicable.

Whenever probing is invoked in \papilo, our \emph{clique probing} procedure is executed first, followed by the standard probing routine. 
In the latter, we exclude variables that have already been probed or fixed by clique probing. Additionally, all bound changes and fixings identified by clique probing are applied locally to the instance used for probing.
Candidate variables are ranked using the same scoring function employed by \papilo{}’s built-in probing routine.

In order to make use of the existing parallelization framework in \papilo, we implemented the option to parallelize the for loops in Line \ref{line::alg::for} in Algorithm \ref{alg::general_algorithm} and Line \ref{line::func::for} in Algorithm \ref{alg::probeclique}.
Each thread receives its own local copy of $B$ and $I$ in order to avoid synchronization overhead for the global data structures $B$ and $I$ in Algorithm~\ref{alg::probeclique}, 
Variable assignments are processed in batches, and after each batch, the local data structures are synchronized sequentially with the global ones.
With this implementation, no locking during the threads is necessary.

In our experiments, however, all algorithms are executed sequentially in order to determine the pure algorithmic advantage of clique probing.
%
Furthermore, \papilo currently has no API to pass implications to an interfaced MIP solver and does not use them internally.
Hence, we do not consider them in our experiments.

\section{Computational Study}
\label{sec::experiments}

In this section, we evaluate how the theoretical speedup translates into practical performance.
In particular, we are interested in investigating the following questions:
\begin{enumerate*}[label=(\alph*)]
\item Is clique probing more effective in the sense that it yields more reductions?
\item How does clique probing affect the overall performance of the solver?
\end{enumerate*}
The first two questions are examined in Section \ref{sec::papilo_experiments} and the latter in Section \ref{sec::scip_experiments}.
First, we explain the experimental setup in Section \ref{sec::setup}.


\subsection{Experimental setup}
\label{sec::setup}

We implemented clique probing in \papilo\cite{GGHpapilo}, which
functions as a standalone presolving library and is also integrated as a default plugin within the presolving phase of the open-source MIP solver \scip.
In our experiments, we run \scip, which invokes \papilo as part of its presolving routine.
It is important to note that \papilo constitutes only one component of the overall presolving process in \scip and may be called multiple times during presolving, with presolving as a whole potentially being called multiple times during the solving process.
Since presolve routines are called before the invocation of \papilo within \scip, \papilo is already executed on a reduced problem and not on the original problem.

\paragraph{Testset}
We base our experiments on the \solver{MIPLIB\,2017} benchmark set~\cite{GleixnerHendelGamrathetal2021MIPLIB2017} consisting of 240~instances and solve each instance with four different random seeds.
We exclude instances that can neither be solved with clique probing enabled nor disabled, on four different seeds, in two hours. 
Further, we exclude instances where \papilo does not identify cliques since they are not affected by the changes.
This results in a total of 93 instances. 
Each seed combination is treated as a separate observation, resulting in a total number of 372 runs.

\paragraph{Software \& Hardware}
For our experiments, we use our implementation in \papilo~3.0 githash \verb|0e1165d6| and compare with the default \papilo githash \verb|32328a5c| as a baseline, which is a verbose variant of \verb|d808e1f3|.
We use a development version of \scip~10.0 with githash \verb|80c1fbfe1b|. 
The experiments were carried out on identical machines with Intel(R) Xeon(R) CPU E7-8880 v4 @ 2.20 GHz.
\paragraph{Parametrization}
By default, \papilo within \scip is executed sequentially, hence we also
used the sequential mode in our experiments.
Since probing large cliques is computationally expensive, we focus only on promising candidate cliques and aim to terminate the probing process early and parameterize Algorithms \ref{alg::general_algorithm} and \ref{alg::probeclique} accordingly.
This approach aligns with the general philosophy of probing.
To avoid long and unsuccessful runtimes, we limit the maximum size of cliques considered to 150 and the maximum number of variables in cliques probed in a single run to 3000. The initial batch size of cliques is set to 2, which determines how many cliques are considered at the start of the process.
Clique probing is aborted if, for two consecutive iterations, fewer than three reductions per propagation are found.
Analogous to probing, we define an additional abortion criterion that terminates clique probing if a threshold of reductions or work limit is exceeded.
If no reduction is found, we also disable clique probing for future rounds of presolving.
Additionally, to avoid probing duplicate variables, we restrict the maximum ratio of probed on 
variables in a clique to 50\%.

\subsection{Analysis of clique probing}
\label{sec::papilo_experiments}
In this section, we analyze the performance of clique probing in terms of the number of fixings, substitutions, and bound changes found in comparison to default \papilo for the first call within \scip. 
A detailed table containing the propagations performed per second, the reductions found per second, and the reductions found for each instance is provided in the repository\footnote{\url{https://github.com/alexhoen/papilo/tree/clique-probing-paper/Paper}}.

\begin{figure}
    \centering

\begin{tikzpicture}
  \begin{loglogaxis}[
    xlabel={Reductions (Default)},
    ylabel={Reductions (Clique)},
    width=\linewidth,
    height=8cm,
    legend style={draw=none},
  ]

    \addplot[
      domain=1:1e6, 
      samples=100,
      black,
      thick,
      dashed
    ] {x};

    \addplot[
      only marks,
      mark=*,
      mark size=1.5pt,
      pbblue,
    ]
    table[row sep=crcr, x=DefFix, y=CliqueFix] {
      DefFix   CliqueFix\\
      3    23\\
0    138\\
0    0\\
0    10\\
24869    24869\\
690    697\\
15297    15206\\
8077    8063\\
191    202\\
143    153\\
23    25\\
1659    1659\\
0    0\\
1919    1919\\
0    0\\
27    168\\
1223    21340\\
0    40\\
1030    237\\
0    0\\
0    0\\
0    71\\
209    477\\
67    361\\
402    387\\
2728    2727\\
18    18\\
0    0\\
1    234\\
0    0\\
    };
    \addplot[
  only marks,
  mark=x,
  mark size=2pt,
  pbgreen,
]
table[row sep=crcr, x=DefBc, y=CliqueBc] {
  DefBc   CliqueBc\\
  0    0\\
0    0\\
0    0\\
0    0\\
80    80\\
8382    8348\\
620    711\\
443    474\\
1    2\\
9    11\\
0    0\\
509    723\\
0    0\\
441    450\\
0    0\\
0    0\\
20    20\\
0    0\\
51    35\\
0    0\\
0    0\\
0    0\\
64    111\\
54    48\\
1196    1173\\
10    10\\
641    641\\
0    0\\
0    0\\
0    0\\
};

\addplot[
  only marks,
  mark=o,
  mark size=2pt,
  pbred,
]
table[row sep=crcr, x=DefSub, y=CliqueSub] {
  DefSub   CliqueSub\\
  0    1\\
0    0\\
0    1\\
0    0\\
0    0\\
0    0\\
0    0\\
0    0\\
0    0\\
0    0\\
0    0\\
2243    2017\\
0    0\\
4026    4026\\
0    0\\
0    0\\
4    4\\
0    0\\
0    0\\
3    9\\
100    101\\
0    30\\
0    26\\
0    2\\
505    514\\
55    55\\
0    0\\
0    41\\
1    3\\
0    1\\
};

  \end{loglogaxis}
\end{tikzpicture}
    \caption{Visualization of reductions of different types, fixings are marked as blue dots, substitutions as red circles, and bound changes as green crosses.}
    \label{fig:papilo_res}
\end{figure}
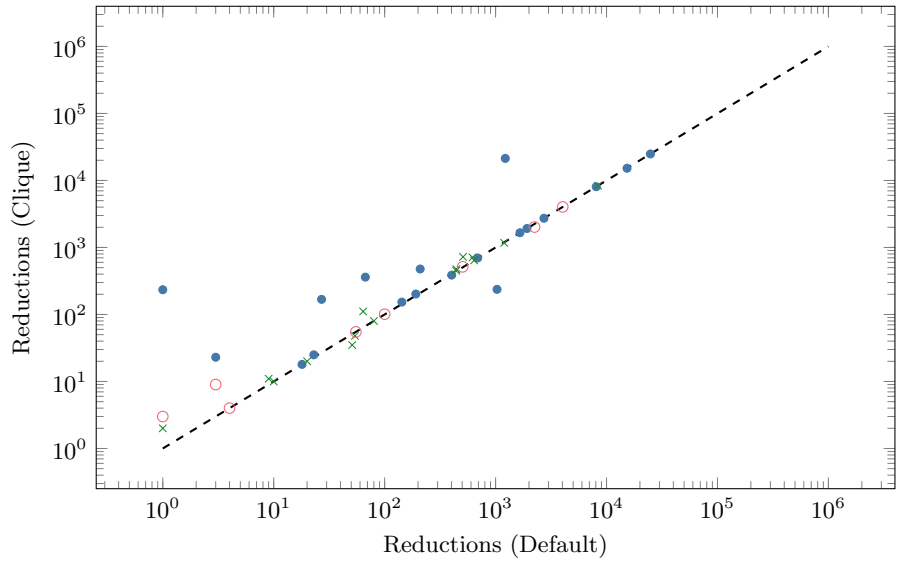

Figure~\ref{fig:papilo_res} visualizes the total number of reductions found, comparing the default \papilo version (horizontal axis) with the clique probing version (vertical axis).
Each data point represents an observation for one instance. Blue points indicate the number of fixings, green crosses show the number of bound changes, and red circles represent substitutions found.
Instances for which the numbers of fixings, substitutions, and bound changes are identical, and the runtime differs by less than 5\% or by at most 0.2 seconds, are omitted from both the table and the figure.
For each of the three categories, the vast majority of data points lie on or above the diagonal, which demonstrates that clique probing is overall more effective than standard probing.

We highlight several instances with particularly pronounced differences.
Two notable negative outliners are the instances \solver{net12} and \solver{square47}: 
On \solver{square47}, the runtime of probing increases from 68.5 to 173.8 seconds.
While default probing aborts probing on the first batch of ten variables without finding any reductions,
Clique probing initially selects two large cliques, which takes a significant amount of time and returns 41 bound changes.
On \solver{net12}, an unfavorable initial clique selection results in premature termination of clique probing, whereas the variable selection of default probing is comparatively ``lucky'': it probes $1600$ variables instead of 250 and finds more reductions. 

In contrast to these outliers, most instances show a clear increase in reductions with only a marginal increase in runtime.
In particular, for \solver{air05}, \solver{neos8}, \solver{neos-860300}, \solver{piperout-08}, \solver{piperout-27}, \solver{supportcase33} 
the improvement is pronounced: the number of fixings increases from 0 to 138, from $1233$ to $21\,340$, from 27 to 168, from 209 to 477, from 67 to 361, and from 1 to 234, respectively, with only a small overhead in runtime.

Generally, the runtime of probing increases slightly when clique probing is active. 
However, the probing process becomes significantly more efficient. 
Compared to using standard probing, \papilo with clique probing finds 5\% more reductions per second of presolving. 
A drawback is that aborting the clique probing process discards all progress made for that clique, which can, for example, on the \solver{square47} instance, increase the overall runtime.


Note that even a small number of additional reductions can substantially impact the behavior of the underlying solver since 
\begin{enumerate*}[label=(\alph*)]
    \item presolvers are executed in a round-based manner to iteratively improve the model and therefore interact with each other, and
    \item contrary to what the name might imply, presolving is invoked multiple times during the overall solution process, for example, when solving sub-MIPs or during restarts.
\end{enumerate*}


\subsection{Impact on overall solver performance}
\label{sec::scip_experiments}
\newcommand{\bracket}[2]{[#1,#2]}

Although the above results are promising, it is difficult to assess solely from the number of reductions and the recorded times how the subsequent solving process behaves with these new changes.
While each reduction achieved through (clique-)probing strengthens the problem formulation, we want to analyze the effectiveness of the reductions obtained by the clique probing approach when integrated within a solver.

In Table \ref{tab::scip_performance}, we compare the performance of clique probing in \scip against the default configuration. All run times are given in seconds and aggregated using the shifted geometric mean with a shift of 1 (sgm1). Nodes are aggregated using the shifted geometric mean with a shift of 100 (sgm100).
The row ``all'' summarizes results over the full test set, while the rows 
$[s,tlim]$ include only instances which were solved by at least one configuration and for which at least one configuration took at least 
$s$ seconds to solve. 
The row ``affected'' lists only those instances where the branch-and-bound tree was altered by clique probing.

Overall, clique probing solves 324 instances compared to 321 in the default setting, thus improving the total number of solved instances by three. We also observe a relative speed-up of approximately 3\% in terms of solving time. 
This improvement becomes more pronounced on harder instances: for instances requiring at least 100 seconds to solve, the speed-up increases to 6\%, and for those requiring at least 1000 seconds, it reaches 8\%. 
The results, therefore, indicate that clique probing is particularly effective on more challenging instances.
\begin{table}
\caption{Performance comparison between \scip with \papilo with clique probing enabled vs. default \papilo. Times are displayed in seconds and aggregated by sgm1. The number of nodes is aggregated by sgm100.}
\label{tab::scip_performance}
\scriptsize
\begin{tabular*}{\textwidth}{@{}l@{\;\;\extracolsep{\fill}}rrrrrrrrr@{}}
\toprule
&           & \multicolumn{3}{c}{clique probing} & \multicolumn{3}{c}{default} & \multicolumn{2}{c}{relative} \\
\cmidrule{3-5} \cmidrule{6-8} \cmidrule{9-10}
Subset                & instances &                                   solved &       time &        nodes & solved&      time &        nodes & time&nodes \\
\midrule
\text{all}&372&324 &310.56 &2125	 & 321 &      318.98	 &2139 & 1.03 &1.01 \\
\text{affected}&  340 & 312 & 350.85 &         3979 & 309 &      372.52 &         4127 &       1.06 &          1.04 \\
\cmidrule{1-10}
\bracket{0}{tlim}& 352 &   324 &  214.30 & 2236 &  321 &  223.45 & 2285 &  1.04 & 1.02 \\
\bracket{1}{tlim}& 352 &   324 &  214.30 & 2236 &  321 &  223.45 & 2285 &  1.04 & 1.02 \\
\bracket{10}{tlim}& 332 &   304 &  330.23	 & 3050 &  301 &  343.43 & 3119 &  1.04 & 1.02 \\
\bracket{100}{tlim}& 264 &   236 &  758.98	 & 7138 &  233 &  804.32 & 7281 &  1.06 & 1.02 \\
\bracket{1000}{tlim}& 176 &   148 &  1987.33	 & 19533 &  145 &  2142.13 & 20793 &  1.08 & 1.06 \\
\bottomrule
\end{tabular*}
\end{table}



\section{Conclusion}
\label{sec::conclusion}

In this paper, we presented clique probing, an approach that exploits the model’s existing data structures to avoid unnecessary propagations during probing while simultaneously providing a stronger relaxation. 
As demonstrated in our experiments, this results in a substantial increase in propagations per second and, consequently, slightly more reductions during probing. 
Overall, this yields a 3\% speed-up in \scip on instances containing cliques. 

There remains, however, room for further improvement.
\begin{enumerate*}[label=(\alph*)]
    \item Since \papilo does not maintain a clique table, and \scip does not share its clique table with \papilo, we are currently limited to detecting only those cliques that appear explicitly in the model; and
    \item there is presently no mechanism for communicating implications back to the underlying solver within \papilo.
\end{enumerate*}

Despite these limitations, the results are encouraging and suggest that similar techniques could be applied to other data structures, such as SOS2 constraints, opening up promising directions for future research.
\\~\\

\textbf{Acknowledgements:}
The work for this article was supported through the Research Campus MODAL funded by the German Federal Ministry of Research, Technology, and Space (fund numbers 05M14ZAM, 05M20ZBM, 05M2025).
Alexander was supported by the Bundesministerium für Bildung und Forschung (16DHBKI071). 

\bigskip


\bibliographystyle{splncs04}
\bibliography{bibliography}
\clearpage
\appendix

\begin{sidewaystable}
\caption{Results of the first probing call of \papilo in \scip. Columns show Default vs. Clique Probing for each metric: Time (s), variable fixings (Fixings), substitutions (Subs), bound changes (BChgs), propagations (Props), average propagations per second (Props/s), and reductions per propagation (Reds/Prop).}
\label{tab:papilo_res}
\resizebox{\textwidth}{!}{
\begin{tabular*}{\textwidth}{@{}l@{\;\;\extracolsep{\fill}} rrrrrrr rrrrrrrr }
    \toprule
    & \multicolumn{2}{c}{Time [s]} 
    & \multicolumn{2}{c}{Fixings} 
    & \multicolumn{2}{c}{Subs} 
    & \multicolumn{2}{c}{BChgs} 
    & \multicolumn{2}{c}{Props}  
    & \multicolumn{2}{c}{Reds/Prop} \\
    \cmidrule(lr){2-3}
    \cmidrule(lr){4-5}
    \cmidrule(lr){6-7}
    \cmidrule(lr){8-9}
    \cmidrule(lr){10-11}
    \cmidrule(lr){12-13}
    Instance 
    & Def & Clique 
    & Def & Clique 
    & Def & Clique 
    & Def & Clique 
    & Def & Clique 
    & Def & Clique  \\
    \midrule
    30n20b8 & 0.2 & 0.2 & 3 & 23 & 0 & 0 & 0 & 1 & 40 & 320 & 0.1 & 0.1 \\
air05 & 0.1 & 0.5 & 0 & 138 & 0 & 0 & 0 & 0 & 20 & 392 & 0.0 & 0.4 \\
bppc4-08 & 0.0 & 0.2 & 0 & 0 & 0 & 0 & 0 & 1 & 30 & 543 & 0.0 & 0.0 \\
brazil3 & 0.2 & 0.3 & 0 & 10 & 0 & 0 & 0 & 0 & 138 & 791 & 0.0 & 0.0 \\
co-100 & 3.2 & 3.9 & 24869 & 24869 & 80 & 80 & 0 & 0 & 62 & 95 & 402.4 & 262.6 \\
cryptanalysiskb128n5obj16 & 5.7 & 6.6 & 690 & 697 & 8382 & 8348 & 0 & 0 & 51316 & 43533 & 0.2 & 0.2 \\
ex10 & 141.5 & 74.2 & 15297 & 15206 & 620 & 711 & 0 & 0 & 32808 & 18831 & 0.5 & 0.8 \\
ex9 & 14.8 & 16.0 & 8077 & 8063 & 443 & 474 & 0 & 0 & 10216 & 9138 & 0.8 & 0.9 \\
mzzv11 & 0.3 & 0.4 & 191 & 202 & 1 & 2 & 0 & 0 & 262 & 398 & 0.7 & 0.5 \\
mzzv42z & 0.2 & 0.3 & 143 & 153 & 9 & 11 & 0 & 0 & 258 & 647 & 0.6 & 0.3 \\
neos-3216931-puriri & 0.1 & 0.2 & 23 & 25 & 0 & 0 & 0 & 0 & 176 & 235 & 0.1 & 0.1 \\
neos-4722843-widden & 4.2 & 10.6 & 1659 & 1659 & 509 & 723 & 2243 & 2017 & 64928 & 70320 & 0.1 & 0.1 \\
neos-5114902-kasavu & 12.2 & 12.9 & 0 & 0 & 0 & 0 & 0 & 0 & 2000 & 2009 & 0.0 & 0.0 \\
neos-5195221-niemur & 0.5 & 1.4 & 1919 & 1919 & 441 & 450 & 4026 & 4026 & 16894 & 19796 & 0.4 & 0.3 \\
neos-662469 & 0.2 & 0.4 & 0 & 0 & 0 & 0 & 0 & 0 & 24 & 299 & 0.0 & 0.0 \\
neos-860300 & 0.8 & 1.5 & 27 & 168 & 0 & 0 & 0 & 0 & 120 & 384 & 0.2 & 0.4 \\
neos8 & 5.9 & 1.9 & 1223 & 21340 & 20 & 20 & 4 & 4 & 3502 & 3459 & 0.4 & 6.2 \\
neos-957323 & 1.1 & 1.9 & 0 & 40 & 0 & 0 & 0 & 0 & 164 & 519 & 0.0 & 0.1 \\
net12 & 0.2 & 0.2 & 1030 & 237 & 51 & 35 & 0 & 0 & 3236 & 921 & 0.3 & 0.3 \\
nexp-150-20-8-5 & 0.1 & 0.1 & 0 & 0 & 0 & 0 & 3 & 9 & 3462 & 3518 & 0.0 & 0.0 \\
ns1644855 & 10.1 & 12.4 & 0 & 0 & 0 & 0 & 100 & 101 & 1950 & 2004 & 0.1 & 0.1 \\
ns1760995 & 6.8 & 31.0 & 0 & 71 & 0 & 0 & 0 & 30 & 266 & 1815 & 0.0 & 0.1 \\
piperout-08 & 0.2 & 0.3 & 209 & 477 & 64 & 111 & 0 & 26 & 682 & 1024 & 0.4 & 0.6 \\
piperout-27 & 0.2 & 0.4 & 67 & 361 & 54 & 48 & 0 & 2 & 266 & 959 & 0.5 & 0.4 \\
rocI-4-11 & 0.2 & 0.2 & 402 & 387 & 1196 & 1173 & 505 & 514 & 6156 & 5733 & 0.3 & 0.4 \\
rocII-5-11 & 2.9 & 2.8 & 2728 & 2727 & 10 & 10 & 55 & 55 & 10068 & 9237 & 0.3 & 0.3 \\
s250r10 & 13.7 & 23.0 & 18 & 18 & 641 & 641 & 0 & 0 & 152 & 420 & 4.3 & 1.6 \\
square47 & 69.5 & 173.8 & 0 & 0 & 0 & 0 & 0 & 41 & 20 & 377 & 0.0 & 0.1 \\
supportcase33 & 0.3 & 1.1 & 1 & 234 & 0 & 0 & 1 & 3 & 62 & 503 & 0.0 & 0.5 \\
wachplan & 0.0 & 0.0 & 0 & 0 & 0 & 0 & 0 & 1 & 22 & 64 & 0.0 & 0.0 \\
    \bottomrule
\end{tabular*}
}
\end{sidewaystable}

\end{document}